\documentclass[11pt, a4paper]{article}
\usepackage{a4wide,amsfonts,amsmath,latexsym,amssymb,euscript}%dsfont,
\usepackage{units,mathrsfs}
\newcommand{\ve}{\varepsilon}
\newcommand{\R}{\mathbb{R}}
\newcommand{\N}{\mathbb{N}}

\newcommand{\fin}
{ \vspace{-0.6cm}
\begin{flushright}
\mbox{$\Box$}
\end{flushright}
\noindent }

\newtheorem{thm}{Theorem}[section]

\newtheorem{prop}[thm]{Proposition}
\newtheorem{lem}[thm]{Lemma}

\newtheorem{rem}[thm]{Remark}
\newtheorem{defi}[thm]{Definition}

\begin{document}

\begin{center}
{\Large{\bf Milstein's type schemes for fractional SDEs}}\\~\\
Mihai Gradinaru\footnote{Institut de Recherche Math{\'e}matique de Rennes, Universit{\'e} de Rennes 1, 
Campus de Beaulieu, 35042 Rennes Cedex, France,
{\tt Mihai.Gradinaru@univ-rennes1.fr }}
and 
Ivan Nourdin\footnote{Laboratoire 
de Probabilit{\'e}s et Mod{\`e}les Al{\'e}atoires, Universit{\'e} Pierre et Marie Curie,
Bo{\^\i}te courrier 188, 4 Place Jussieu, 75252 Paris Cedex 5, France,
{\tt ivan.nourdin@upmc.fr}}
\end{center}

{\small
\noindent
{\bf R\'esum\'e.} On \'etudie la vitesse exacte de convergence de certains sch\'emas d'approximation associ\'es \`a des \'equations 
diff\'erentielles stochastiques scalaires dirig\'ees par le mouvement brownien fractionnaire $B$. On utilise le comportement 
asymptotique des variations \`a poids de $B$, et la limite de l'erreur entre la solution et son approximation est calcul\'ee 
de fa\c{c}on explicite.\\

\noindent
{\bf Abstract: }
Weighted  power variations of fractional Brownian motion $B$ are used to compute the exact rate of convergence 
of some approximating schemes associated to one-dimensional stochastic differential equations (SDEs) driven by $B$. The limit 
of the error between the exact solution and the considered scheme 
is computed explicitly.\\


}

\section{Introduction}
\setcounter{footnote}{0}

Let $B=(B_t)_{t\in [0,1]}$ be a fractional Brownian motion with Hurst index $H\in (0,1)$.
That is, $B$ is a centered Gaussian process with covariance function given by
$$
{\rm Cov}(B_{s},B_{t})=\frac{1}{2}(s^{2H}+t^{2H}-|t-s|^{2H}),\quad s,t\in [0,1].
$$
For $H=\nicefrac{1}{2}$, $B$ is a standard Brownian motion, while for $H\neq\nicefrac{1}{2}$, it 
is neither a semimartingale, nor a Markov process. Moreover, it holds, for any $p>1$:
$$
{\rm E}|B_{t}-B_{s}|^{p}=c_{p}|t-s|^{pH},\quad s,t\in[0,1],\quad\mbox{with }
\quad c_{p}={\rm E}(|G|^{p}),\quad G\sim\mathscr{N}(0,1),
$$
and, consequently, almost all sample paths of $B$ are H{\"o}lder continuous of any order 
$\alpha\in (0,H)$. 

The study of stochastic differential equations driven by $B$
has been considered by using several methods. 
For instance, in \cite{NR} one uses fractional calculus of same type as in \cite{zahle}; in \cite{CQ} one uses 
rough paths theory introduced in \cite{lyons},
and in \cite{N06} one uses 
regularization method used firstly in  
\cite{RV93}.

In the present paper, we consider the easiest stochastic differential equation
involving fractional Brownian motion, that is
\begin{equation}\label{eqintro}
dX_t=\sigma(X_t)dB_t
,\,t\in [0,1],\quad X_0=x\in\R.
\end{equation} 
Here and in the rest of the paper, $\sigma\in{\rm C}^{\infty}(\R)$ 
stands for a real function which is bounded with bounded derivatives. Let us denote by 
$\phi:\R^2\rightarrow\R$ the flow associated to $\sigma$, that is the unique solution to 
\begin{equation}\label{flow}
\phi(x,y)=x+\int_{0}^{y}\sigma(\phi(x,z))dz,\quad x,y\in\R.
\end{equation}
Assume that the integral with respect to 
$B$ we consider in (\ref{eqintro}) verifies the following 
It{\^o}-Stratonovich type formula:
\begin{equation}\label{dol}
f(B_t)=f(0)+\int_0^t f'(B_s)dB_s,\,\quad t\in[0,1],\quad
f:\mathbb{R}\to\mathbb{R}\mbox{ smooth enough.}
\end{equation}
Then, combined with (\ref{flow}), one easily checks that
\begin{equation}\label{X}
X_t^x=\phi(x,B_t),\quad t\in[0,1],
\end{equation}
is a solution to (\ref{eqintro}).

Approximating schemes for stochastic differential equations of the type 
\begin{equation}\label{deriveqintro}
dX_t=\sigma(X_t)dB_t+b(X_t)dt
,\,t\in [0,1],\quad X_0=x\in\R,
\end{equation} 
have been considered only in few articles. The first work in that direction is 
\cite{lin}. Precisely, whenever $H>\nicefrac{1}{2}$, it is shown that the Euler approximation of equation (\ref{deriveqintro}) -- but in the particular case where 
$\sigma(X_t)$ is replaced by $\sigma(t)$, that is the so-called additive case --
converges uniformly in probability. In \cite{N05} one introduces (see also \cite{talay}) 
some approximating schemes for the analogue of (\ref{deriveqintro}) where $B$ is replaced by 
any H{\"o}lder continuous function. One determines 
upper error bounds and, in particular, these results apply almost surely when the driving 
H{\"o}lder continuous function is a single path of the fractional Brownian motion $B$, and this 
for any Hurst index $H\in(0,1)$. In \cite{mishura}, upper error bounds for Euler approximations of solutions of (\ref{deriveqintro}) are derived whenever $H$ is bigger than $\nicefrac12$. The convergence of Euler schemes has also been studied in \cite{davie} in the context of the rough paths theory.

Results on lower error bounds are available only since very recently: see \cite{neuenkirch} 
for the additive case, and \cite{NN} for equation 
(\ref{deriveqintro}) (see also \cite{ANopt} where approximation methods
with respect to a mean square error are analysed). More precisely, it is proved
in \cite{NN} that the Euler scheme $\widetilde{X}
=\{\widetilde{X}^{(n)}\}_{n\in\N}$ associated to (\ref{deriveqintro}) verifies, 
under some classical assumptions on $\sigma$ and $b$ 
and whenever $H\in(\nicefrac{1}{2},1)$, that
\begin{equation}\label{euler-nn}
n^{2H-1}\left[ \widetilde{X}^{(n)}_1 - X_1 \right]  
\,{\stackrel{{\rm a.s.}}{\longrightarrow}}\, 
-\frac{1}{2}\int_0^1 \sigma'(X_s)D_sX_1ds,\quad\mbox{as $n\to\infty$}.
\end{equation}
Here, $D_sX_1$ denotes the Malliavin derivatives of $X_1$ with respect to $B$.
Observe that the upper and lower error bounds are obtained from 
an almost sure convergence, which is somewhat surprising when compared with the 
case $H=\nicefrac{1}{2}$, see below.
In \cite{NN}, it is proved that, for the so-called Crank-Nicholson scheme 
$\overline{X}=\{\overline{X}^{(n)}\}_{n\in\N}$ associated to (\ref{eqintro}) 
and defined by
\begin{equation}\label{crank-nn}
\left\{
\begin{array}{lll}
\overline{X}^n_0=x,\\
\overline{X}^{(n)}_{(\ell+1)/n}=\overline{X}^{(n)}_{\ell/n}+\frac{1}{2}\left(\sigma(\overline{X}^{(n)}_{\ell/n})
+\sigma(\overline{X}^{(n)}_{(\ell+1)/n})\right)(B_{(\ell+1)/n}-B_{\ell/n}),\quad 
\ell\in\{0,\ldots,n-1\},
\end{array}
\right.
\end{equation}
the following convergence holds
for $\sigma$ regular enough and whenever $H\in(\nicefrac{1}{3},\nicefrac{1}{2})$:
\begin{equation}\label{nn1}
n^{\alpha}\left[ 
\overline{X}^{(n)}_1 
- X_1 \right]  
\,{\stackrel{{\rm Prob}}{\longrightarrow}}\, 
0,\quad\forall\alpha<3H-\nicefrac{1}{2}, \quad\mbox{as $n\to\infty$}.
\end{equation}
In the particular case where the diffusion coefficient $\sigma$ verifies
$\sigma(x)^2= \alpha x^2 + \beta x+\gamma$, for some $\alpha,\beta,\gamma\in\R$, one can derive the
exact rate of convergence and one proves that, as $n\to\infty$:
\begin{equation}\label{nn2}
n^{3H-\nicefrac{1}{2}}\left[ \overline{X}^{(n)}_1 - X_1 \right]  
\,{\stackrel{{\rm Law}}{\longrightarrow}}\, 
\frac{\alpha}{12}\sigma(X_1)\,G.
\end{equation}
Here, $G$ is a centered Gaussian random variable independent of $X_1$, whose variance depends uniquely on $H$. In particular, the upper and lower error bounds are obtained here from a convergence in law.

As we said, the convergence in (\ref{euler-nn}) is  
somewhat surprising, since there is no analogue for the case of the standard Brownian motion. 
More precisely, when $H=\nicefrac{1}{2}$, it is proved in \cite{KP} that
the Euler scheme (\ref{euler-nn}) verifies (by denoting $X^\textrm{It{\^o}}$ the solution
of (\ref{deriveqintro}) in the It{\^o} sense), 
\begin{equation}\label{euler1/2}
\sqrt{n}\left[ \widetilde{X}^{(n)}_1 - X_1^{\textrm{It{\^o}}} \right]  
\,{\stackrel{{\rm Law}}{\longrightarrow}}\, 
-\frac{1}{\sqrt{2}}\,Y_1\,\int_0^1 \sigma(X_{s}^{\textrm{Ito}})\sigma'(X_s^{\textrm{Ito}})Y_s^{-1}dW_s,\quad\mbox{as $n\to\infty$}.
\end{equation}
Here, $W$ is a Brownian motion independent of the Brownian motion $B$ and
$$ 
Y_{t}= \exp \left(
\int_{0}^{t} b'(X_{s}^{\textrm{Ito}}) ds
-\frac{1}{2}\int_{0}^{t} \sigma\sigma'(X_{s}^{\textrm{Ito}})\, ds 
+ \int_{0}^{t} \sigma'(X_{s}^{\textrm{Ito}})\, dB_{s}  \right), \quad t \in [0,1].$$ 
On the other hand, it can be proved (see Remark \ref{rm34}.2) that,
for the Crank-Nicholson scheme (\ref{crank-nn}), we have, as $n\to\infty$:
\begin{equation}\label{nn3}
n\left[ \overline{X}^{(n)}_1 - X_1^{\textrm{Str}}  \right] 
\,{\stackrel{{\rm Law}}{\longrightarrow}}\, 
\frac{1}{24}\,\int_0^1 (\sigma^2)''(X^{\textrm{Str}}_s)[\sqrt{15}dW_s
+3dB_s]
+\frac{1}{16}\,\int_0^1 \sigma(\sigma^2)'''(X^{\textrm{Str}}_s)ds,
\end{equation}
where $X^\textrm{Str}$ denotes the solution of (\ref{eqintro}) in the Stratonovich sense.

In the present paper, we are interested in a better understanding of the phenomenoms
observed in (\ref{euler-nn}), (\ref{nn2}), (\ref{euler1/2}) or (\ref{nn3}).
What type of convergence allows to derive the upper and lower error bounds
for some natural scheme of Milstein's type?
More precisely, let us define, by induction, the family of differential operators 
$(\mathscr{D}^j)_{j\in\N\cup\{0\}}$ as
\begin{equation}\label{do}
\mathscr{D}^0f=f,\quad\mathscr{D}^1f=f'\sigma\quad\mbox{ and, for $j\ge 2$, }
\quad\mathscr{D}^jf=\mathscr{D}^1(\mathscr{D}^{j-1}f).
\end{equation}
For instance, the first $\mathscr{D}^j\sigma$'s are given by:
$$
\mathscr{D}^0\sigma=\sigma,\,
\mathscr{D}^1\sigma=\sigma\sigma',\,
\mathscr{D}^2\sigma=\sigma\sigma'^2+\sigma^2\sigma'',\,
\mathscr{D}^3\sigma=\sigma\sigma'^3+4\sigma^2\sigma'\sigma''+\sigma^3\sigma''',\,
\mbox{etc.}
$$
Now, let us consider the following scheme introduced in \cite{N05}:
\begin{equation}\label{scheme}
\left\{
\begin{array}{lll}
\widehat{X}^{(n)}_{0}=x\\
\widehat{X}^{(n)}_{\nicefrac{(\ell+1)}{n}}=\widehat{X}^{(n)}_{\nicefrac{\ell}{n}}+
\sum_{j=0}^m\,\frac{1}{(j+1)!}
\,\mathscr{D}^j\sigma(\widehat{X}^{(n)}_{\nicefrac{\ell}{n}})
\,\big(\Delta B_{\nicefrac{\ell}{n}})^{j+1}
,\quad \ell\in\{0,\ldots,n-1\},
\end{array}
\right.
\end{equation}
the integer $m\in\N\cup\{0\}$ being called the size of 
$\widehat{X}=\{\widehat{X}^{(n)}\}_{n\in\N}$. 
Here, for $j,n\in\N$ and $\ell\in\{0,\ldots,n-1\}$, we set 
$\Delta B_{\nicefrac{\ell}{n}}$ instead of $B_{\nicefrac{(\ell+1)}{n}}-B_{\nicefrac{\ell}{n}}$ for simplicity. 
The idea of introducing these schemes will be better explained in Section 3 below.
For the moment, just observe that Euler (resp. Milstein) scheme corresponds to
$m=0$ (resp. $m=1$).

The aim of the present paper is to answer the following questions. Does the sequence 
$\{\widehat{X}^{(n)}_1\}_{n\in\N}$ converge? Is the limit $X^x_1$ given by (\ref{X}), as could be reasonably 
expected? 
What is the rate of convergence? Are upper and lower error bounds
obtained from a convergence in law or rather from a pathwise type convergence?

The paper is organized as follows: the next section reviews 
some very recent results concerning 
the asymptotic behavior of weighted power variations of fractional Brownian motion. In 
section 3,  after recalling the definition and the main properties of 
the so-called Newton-Cotes 
integral, we explain how to use it in order to study (\ref{eqintro}).
Finally, in section 4, we state and prove our results concerning the exact rate of 
convergence associated to (\ref{scheme}).

\section{Asymptotic behavior of weighted power variations}
\setcounter{equation}{0}

Let $\kappa\geq 2$ be an integer, and $h,g:\R\to\R$ be two functions belonging
to $\mathrm{C}^\infty$. Assume moreover that $h$ and $g$ are bounded with bounded derivatives.
Denote by $\mu_{2n}$ the $2n$-moment of a random variable 
$G\sim\mathscr{N}(0,1)$. The following theorem collects some very recent results about the asymptotic behavior
of the so-called {\it weighted power variations} of $B$, defined by 
\begin{equation*}
\sum_{\ell=0}^{n-1} h(B_{\nicefrac{\ell}{n}}) \big(
\Delta B_{\nicefrac{\ell}{n}}
\big)^\kappa
\quad\mbox{(recall that $\Delta B_{\nicefrac{\ell}{n}}$ stands for $B_{\nicefrac{(\ell+1)}{n}}-B_{\nicefrac{\ell}{n}}$).}
\end{equation*}

\begin{thm}\label{pwfbm}
\begin{enumerate}
\item If $\kappa$ is even and $H\in (0,1)$ then, as $n\to\infty$:
\begin{equation}\label{2-pow}
n^{\kappa H-1}\sum_{\ell=0}^{n-1} h(B_{\nicefrac{\ell}{n}})
(\Delta B_{\nicefrac{\ell}{n}})^{\kappa }
\,\,\,{\stackrel{{\rm Prob}}{\longrightarrow}}\,\,\,
\mu_\kappa
\int_0^1 h(B_s) ds.
\end{equation}
\item If $\kappa$ is odd and $H\in (0,\nicefrac{1}{2})$ then, as 
$n\to\infty$:
\begin{equation}\label{3-pow}
n^{(\kappa+1) H-1}
\sum_{\ell=0}^{n-1} h(B_{\nicefrac{\ell}{n}})(\Delta B_{\nicefrac{\ell}{n}})^{\kappa }
\,\,\,{\stackrel{{\rm Prob}}{\longrightarrow}}\,\,\,
-\frac{\mu_{\kappa+1}}{2}
\int_0^1 h'(B_s) ds.
\end{equation}
\item If $\kappa$ is odd and $H=\nicefrac{1}{2}$ then, as $n\to\infty$,
\begin{multline}\label{mb3-pow}
\left(B_1,
n^{\frac{\kappa -1}{2}}\sum_{\ell=0}^{n-1} 
g(B_{\nicefrac{\ell}{n}})(\Delta B_{\nicefrac{\ell}{n}})^{\kappa },
n^{\frac{\kappa -1}{2}}\sum_{\ell=0}^{n-1} 
h(B_{\nicefrac{\ell}{n}})(\Delta B_{\nicefrac{\ell}{n}})^{\kappa+1 }
\right)\\
{\stackrel{{\rm Law}}{\longrightarrow}}\,\,\,
\left(B_1,
\int_0^1 g(B_s) \big(\sqrt{\mu_{2\kappa}}dW_s
+\mu_{\kappa+1}dB_s\big)
,
\mu_{\kappa+1}
\int_0^1 h(B_s) ds
\right),
\end{multline}
with $W$ another standard Brownian motion independent of $B$.
\item If $\kappa$ is odd and $H\in (\nicefrac{1}{2},1)$ then,
as $n\to\infty$:
\begin{equation}\label{28bis}
n^{(\kappa-1) H}
\sum_{\ell=0}^{n-1} h(B_{\nicefrac{\ell}{n}})(\Delta B_{\nicefrac{\ell}{n}})^{\kappa }
\,\,\,{\stackrel{{\rm Prob}}{\longrightarrow}}\,\,\,
\mu_{\kappa+1}
\int_0^{B_1} h(x) ds.
\end{equation}
\end{enumerate}
\end{thm}
\begin{rem}
{\rm 
\begin{enumerate}
\item For sake of conciseness, we omit the proof of Theorem \ref{pwfbm}. 
We give below some ideas and references for the proofs.
\item The convergence (\ref{2-pow}) is actually almost sure. 
Its proof is a classical result when $h\equiv 1$ (see e.g. \cite{KG} 
when $\kappa=2$). If $h$ is arbitrary, the proof could be completed along the lines of the proof of 
Lemma 3.1, p. 7-8 in \cite{GN03}. 
\item Proofs of (\ref{3-pow}) and (\ref{28bis}) can be completed along the lines of \cite[Corollary 2]{NNT}.
\item We shall see that, for the standard Brownian motion case, in order to study our Milstein's type schemes 
one needs the behaviour of the triplet (\ref{mb3-pow}) and not only the behaviour of 
the second coordinate in (\ref{mb3-pow}).
The proof of (\ref{mb3-pow}) can be completed along the lines of \cite[Corollary 2.9]{NoPe}. More precisely, 
using the methodology introduced in this latter reference, we first prove that\footnote{
Using the notion of {\it stable convergence} 
for random variables, (\ref{oil}) is equivalent to say that
$$
n^{\frac{\kappa -1}{2}}\sum_{\ell=0}^{n-1} 
g(B_{\nicefrac{\ell}{n}})(\Delta B_{\nicefrac{\ell}{n}})^{\kappa }
\quad{\stackrel{\mathscr{F}^{B}-{\rm stably}}{\longrightarrow}}
\int_0^1 g(B_s) \big(\sqrt{\mu_{2\kappa}}dW_s
+\mu_{\kappa+1}dB_s\big).
$$
Here, $\mathscr{F}^{B}$ denotes the $\sigma$-field generated by $(B_t)_{t\in[0,1]}$
(see also Theorem 1.1, p. 3 in \cite{jacod}).
}, as $n\to\infty$,
\begin{multline}\label{oil}
\left((B_t)_{t\in[0,1]},n^{\frac{\kappa -1}{2}}\sum_{\ell=0}^{n-1} 
g(B_{\nicefrac{\ell}{n}})(\Delta B_{\nicefrac{\ell}{n}})^{\kappa }\right)\\
{\stackrel{{\rm Law}}{\longrightarrow}}
\left((B_t)_{t\in[0,1]},
\int_0^1 g(B_s) \big(\sqrt{\mu_{2\kappa}}dW_s
+\mu_{\kappa+1}dB_s\big)\right).
\end{multline}
Then, using the fact that  (see (\ref{2-pow})), as as $n\to\infty$,
$$
n^{\frac{\kappa -1}{2}}\sum_{\ell=0}^{n-1} 
h(B_{\nicefrac{\ell}{n}})(\Delta B_{\nicefrac{\ell}{n}})^{\kappa+1 }
{\stackrel{{\rm Prob}}{\longrightarrow}}\,\,\,
\mu_{\kappa+1}
\int_0^1 h(B_s) ds,
$$
and that $\mu_{\kappa+1}
\int_0^1 h(B_s) ds$ is a random variable measurable with respect to $B$, the desired
conclusion follows easily.
\item Other results on weighted variations of fractional Brownian motion (or related processes) 
can be found in \cite{CNW} and \cite{LL}.
\end{enumerate}
}
\end{rem}

\section{Newton-Cotes integral and fractional SDEs}
\setcounter{equation}{0}

In the sequel, we will use, as integral with respect to $B$, the so-called Newton-Cotes 
integral introduced in \cite{GNRV} and studied further in \cite{N06}:  
\begin{defi}
Let $f:\mathbb{R}\to\mathbb{R}$ be a continuous function, $X,Y$ be two continuous processes 
on $[0,1]$ and $N\in\mathbb{N}\cup\{0\}$. The $N$-{\sl order Newton-Cotes integral of} $f(Y)$ {\sl with 
respect to} $X$ is defined by:
\begin{equation}\label{nc-def}
\int_{0}^{t}f(Y_{s})d^{{\rm NC},N}X_{s}
:=\lim_{\ve\downarrow 0} {\rm prob}\frac{1}{\ve}
\int_{0}^{t}ds (X_{s+\ve}-X_{s})
\int_{0}^{1}f(Y_{s}+\alpha(Y_{s+\ve}-Y_{s}))\nu_{N}(d\alpha),\,\,t\in[0,1],
\end{equation}
provided the limit exists. Here $\nu_{0}=\delta_{0}$, $\nu_{1}=\nicefrac{(\delta_{0}+\delta_{1})}{2}$ 
and, for $N\geq 2$, 
\begin{equation*}
\nu_{N}=\sum_{j=0}^{2N-2}\left(\int_{0}^{1}\prod_{k\neq j}
\frac{2(N-1)u-k}{j-k}\,du\right)\delta_{\nicefrac{j}{(2N-2)}},
\end{equation*}
$\delta_{a}$ being the Dirac measure at point $a$. 
\end{defi}

\begin{rem}
{\rm  
\begin{enumerate}
\item The 0- and 1-order Newton-Cotes integrals
are nothing but the forward integral and the symmetric integral  
in the sense of Russo-Vallois \cite{RV93}, respectively:
\begin{equation*}
\int_{0}^{t}f(Y_{s})d^{{\rm NC},0}X_{s}=
\int_{0}^{t}f(Y_{s})d^{-}X_{s}=
\lim_{\ve\downarrow 0}{\rm prob}\frac{1}{\ve}
\int_{0}^{t}f(Y_{s})(X_{s+\ve}-X_{s})ds,
\end{equation*}
and 
\begin{equation*}
\int_{0}^{t}f(Y_{s})d^{{\rm NC},1}X_{s}=
\int_{0}^{t}f(Y_{s})d^{\circ}X_{s}=
\lim_{\ve\downarrow 0}{\rm prob}\frac{1}{\ve}
\int_{0}^{t}\frac{f(Y_{s+\ve})+f(Y_{s})}{2}(X_{s+\ve}-X_{s})ds.
\end{equation*}
\item Another way to define $\nu_N$ is to view it as the unique discrete signed probability 
carried by $\nicefrac{j}{(2N-2)}$ ($j=0,\ldots,2N-2$), which coincides with Lebesgue 
measure on polynomials of degree smaller than $2N-1$.
\end{enumerate}
}
\end{rem}

The Newton-Cotes integral defined by (\ref{nc-def}) is actually a special case of so-called 
$N$-order $\nu$-integrals introduced in \cite{GNRV}, p. 789. Moreover, in the same cited paper, 
p. 795, one proves that the $N$-order 
Newton-Cotes integral of $f(B)$ with respect to $B$ exists for any 
$f\in{\rm C}^{4N+1}$ if and only if $H\in(\nicefrac{1}{(4N+2)},1)$. In this case, an It{\^o}'s type 
change of variables formula holds: for 
any antiderivative $F$ of $f$, we can write 
\begin{equation}\label{ito-m}
F(B_{t})-F(0)=\int_{0}^{t}f(B_{s})d^{{\rm NC},N}B_{s},\quad t\in[0,1].
\end{equation}
Moreover, as a consequence of (\ref{ito-m}), let us note that 
$$
\int_0^t f(B_s)d^{{\rm NC},N} B_s
=\int_0^t f(B_s)d^{{\rm NC},n} B_s=F(B_t)-F(0),
$$ 
as soon as  $f\in{\rm C}^{4N+1}$, $n<N$ and $H\in (\nicefrac{1}{(4n+2)},1)$. Therefore, for 
$f$ regular enough, it is possible to define the Newton-Cotes integral 
without ambiguity by:
\begin{equation}\label{nc-gen}
\int_0^t f(B_s)d^{{\rm NC}} B_s:=
\int_0^t f(B_s)d^{{\rm NC},n} B_s\quad\quad\mbox{if }H\in (\nicefrac{1}{(4n+2)},1).
\end{equation}
Set $n_{_{H}}:=\inf\{n\ge 1:\,H>\nicefrac{1}{(4n+2)}\}$.
Hence, an immediate consequence of (\ref{ito-m}) and (\ref{nc-gen}) is that, for any $H\in (0,1)$
and any $f:\R\rightarrow\R$ of class ${\rm C}^{4n_{_{H}}+1}$, the following It{\^o}'s type 
change of variables formula holds:
\begin{equation}\label{ito}
F(B_t)=F(0)+\int_0^t f(B_s)d^{{\rm NC}} B_s,
\quad\mbox{ for any antiderivative $F$ of $f$.}
\end{equation}

\begin{rem}
{\rm 
In the sequel 
we will only use the fact that the Newton-Cotes integral verifies the 
classical change of variable formula (\ref{dol}) or (\ref{ito}). Consequently, any other stochastic integral verifying
(\ref{ito}) could be used in the following.
}
\end{rem}

All along this paper we will work with an ellipticity assumption, and we will
also need regularity for the function 
$\sigma$. More precisely, we suppose
$$
({\mathscr E})\quad\quad\mbox{$\inf_\mathbb{R}|\sigma|>0$}\quad
\mbox{and $\sigma\in{\rm C}^{\infty}(\mathbb{R})$ is bounded with bounded derivatives}.
$$
Under hypothesis $({\mathscr E})$,  the flow $\phi$ associated to $\sigma$, given by 
(\ref{flow}), is well-defined and verifies the group property: 
\begin{equation}\label{sgp}
\forall x,y,z\in\mathbb{R},\quad
\phi(\phi(x,y),z)=\phi(x,y+z).
\end{equation}
Note that the process $X^x$ given by (\ref{X}) verifies:
\begin{equation}\label{fsde-int}
X^x_{t}=x+\int_0^t \sigma(X^x_{s})d^{{\rm NC}}B_{s},\quad t\in [0,1], 
\end{equation}
as we can see immediately, by applying (\ref{ito}).
\begin{rem}
{\rm
In \cite{N06} (see also \cite{NS}), one studies a notion of solution for 
(\ref{fsde-int}) and also the existence and the uniqueness of solution. Note however that, in 
the present work, we will only use the fact that there exists a natural solution to (\ref{fsde-int}) 
given by (\ref{X}).
}
\end{rem}

The following result explains the definition (\ref{scheme}). By using (\ref{do}), the process $X^x$ 
defined by (\ref{X}) can be expanded as follows:
\begin{lem}\label{lm1}
For any integers $m\ge 0$, $n\ge 1$ and $\ell\in\{0,\ldots,n-1\}$, 
we have
\begin{multline}\label{dvpx}
X^x_{\nicefrac{(\ell+1)}{n}}
=
X^x_{\nicefrac{\ell}{n}}+\sum_{j=0}^m\,\frac{1}{(j+1)!}
\,\mathscr{D}^j\sigma(X^x_{\nicefrac{\ell}{n}})
\,\big(\Delta B_{\nicefrac{\ell}{n}}\big)^{j+1}\\
+
\int_{\nicefrac{\ell}{n}}^{\nicefrac{(\ell+1)}{n}} d^{{\rm NC}}B_{t_1}
\int_{\nicefrac{\ell}{n}}^{t_1} d^{{\rm NC}}B_{t_2}
\ldots
\int_{\nicefrac{\ell}{n}}^{t_m} d^{{\rm NC}}B_{t_{m+1}}
\int_{\nicefrac{\ell}{n}}^{t_{m+1}} \mathscr{D}^{m+1}\sigma(X^x_{t_{m+2}}) d^{{\rm NC}}B_{t_{m+2}}.\\
\end{multline}
\end{lem}

\noindent
{\bf Proof.}
We proceed by induction on $m$.
By applying (\ref{fsde-int}), and then using (\ref{X}) and (\ref{ito}), we can write:
\begin{multline*}
X^x_{\nicefrac{(\ell+1)}{n}}
=X^x_{\nicefrac{\ell}{n}}+
\sigma(X^x_{\nicefrac{\ell}{n}})
\Delta B_{\nicefrac{\ell}{n}}
+\int^{{\nicefrac{(\ell+1)}{n}}}_{{\nicefrac{\ell}{n}}} (\sigma(X^x_{t_1})-\sigma(X^x_{\nicefrac{\ell}{n}})\big)
d^{{\rm NC}}B_{t_1}\\
=X^x_{\nicefrac{\ell}{n}}+
\sigma(X^x_{\nicefrac{\ell}{n}})
\Delta B_{\nicefrac{\ell}{n}}
+\int^{{\nicefrac{(\ell+1)}{n}}}_{{\nicefrac{\ell}{n}}} d^{{\rm NC}}B_{t_1}
\int^{t_1}_{{\nicefrac{\ell}{n}}} \sigma\sigma'(X^x_{t_2})
d^{{\rm NC}}B_{t_2},
\end{multline*}
which is exactly (\ref{dvpx}) for $m=0$.

Now, let us assume that (\ref{dvpx}) is true for some 
$m\in\mathbb{N}\cup\{0\}$. Then we can write
\begin{multline}\label{eq79}
X^x_{\nicefrac{(\ell+1)}{n}}=X^x_{\nicefrac{\ell}{n}}+\sum_{j=0}^m\,\frac{1}{(j+1)!}
\,\mathscr{D}^j\sigma(X^x_{\nicefrac{\ell}{n}})
\,\big(\Delta B_{\nicefrac{\ell}{n}})^{j+1}\\
+
\mathscr{D}^{m+1}\sigma(X^x_{\nicefrac{\ell}{n}}) 
\int_{\nicefrac{\ell}{n}}^{\nicefrac{(\ell+1)}{n}} d^{{\rm NC}}B_{t_1}
\ldots
\int_{\nicefrac{\ell}{n}}^{t_m} d^{{\rm NC}}B_{t_{m+1}}
\int_{\nicefrac{\ell}{n}}^{t_{m+1}} 
d^{{\rm NC}}B_{t_{m+2}}\\
+
\int_{\nicefrac{\ell}{n}}^{\nicefrac{(\ell+1)}{n}} d^{{\rm NC}}B_{t_1}
\ldots
\int_{\nicefrac{\ell}{n}}^{t_m} d^{{\rm NC}}B_{t_{m+1}}
\int_{\nicefrac{\ell}{n}}^{t_{m+1}} 
\big(\mathscr{D}^{m+1}
\sigma(X^x_{t_{m+2}}) 
-
\mathscr{D}^{m+1}
\sigma(X^x_{\nicefrac{\ell}{n}})\big)
d^{{\rm NC}}B_{t_{m+2}}.
\end{multline}
On one hand, using (\ref{ito}) repeatedly, it is immediate to compute that
$$
\int_{\nicefrac{\ell}{n}}^{\nicefrac{(\ell+1)}{n}} d^{{\rm NC}}B_{t_1}
\ldots
\int_{\nicefrac{\ell}{n}}^{t_m} d^{{\rm NC}}B_{t_{m+1}}
\int_{\nicefrac{\ell}{n}}^{t_{m+1}} 
d^{{\rm NC}}B_{t_{m+2}}=\frac{1}{(m+2)!}\big(\Delta B_{\nicefrac{\ell}{n}}\big)^{m+2}.
$$ 
On the other hand, using (\ref{X}) and again (\ref{ito}), we can write
$$
\mathscr{D}^{m+1}
\sigma(X^x_{t_{m+2}}) 
-
\mathscr{D}^{m+1}
\sigma(X^x_{\nicefrac{\ell}{n}})
=
\int_{\nicefrac{\ell}{n}}^{t_{m+2}} 
\sigma\big(\mathscr{D}^{m+1}\sigma\big)'
(X^x_{t_{m+3}}) d^{{\rm NC}}B_{t_{m+3}}.
$$ 
Finally, putting these latter two equalities in (\ref{eq79}) and noting that
$
\sigma\big(\mathscr{D}^{m+1}\sigma\big)'
=\mathscr{D}^{m+2}\sigma$ by definition, we obtain that (\ref{dvpx}) is true
for $m+1$. The proof by induction is done. \fin

Clearly, (\ref{scheme}) is the natural scheme constructed from (\ref{dvpx}), by considering
the $(m+2)$th multiple integral in the right hand side of (\ref{dvpx}) as a remainder.

\section{Rate of convergence of the approximating schemes}
\setcounter{equation}{0}

\subsection{Statement of the main result}

Recall that we denote by $\mu_{2n}$ the $2n$-moment of a random variable 
$G\sim\mathscr{N}(0,1)$.
For $m\in\N\cup\{0\}$, let us introduce the functions $g_m,h_m:\R\to\R$ given by:
\begin{equation}\label{gm}
g_m=-\sigma'\,h_m+h_{m+1}
\quad\mbox{and}\quad
h_m=
-\frac{\nicefrac{(\mathscr{D}^{m+1}\sigma)}{\sigma}}{(m+2)!}.
\end{equation}

Our main result contains a complete answer to the questions in the 
introduction and can be stated as follows:
\begin{thm}\label{main-result}
Assume that hypothesis $({\mathscr E})$ is in order, and let $m\in\mathbb{N}\cup\{0\}$. Then, for any $H\in(\nicefrac{1}{(m+2)},1)$, the sequence $\{\widehat{X}^{(n)}_1\}_{n\in\mathbb{N}}$ 
defined by (\ref{scheme}) converges almost surely toward $X^x_1=\phi(x,B_1)$
as $n\to\infty$.
Moreover,

\textbullet\quad when $m$ is even and $H\in(\nicefrac{1}{(m+2)},1)$, 
\begin{equation}\label{thm1-eq1}
n^{(m+2)H-1}\left[\widehat{X}^{(n)}_1-X^x_1\right]
\,\,\,{\stackrel{{\rm Prob}}{\longrightarrow}}\,\,\,
\mu_{m+2}\,\sigma(X^x_1)\int_0^1 h_{m}(X_s^x)ds\,;
\end{equation}

\textbullet\quad when $m$ is odd and  $H\in(\nicefrac{1}{(m+2)},\nicefrac{1}{2})$,
\begin{equation}\label{thm1-eq2}
n^{(m+3)H-1}\left[\widehat{X}^{(n)}_1-X^x_1\right]
\,\,\,{\stackrel{{\rm Prob}}{\longrightarrow}}\,\,\,
\mu_{m+3}\,\sigma(X^x_1)\int_0^1 \big(g_{m}-\frac12\sigma h'_m\big)(X_s^x)ds\,;
\end{equation}

\textbullet\quad when $m$ is odd and  $H=\nicefrac{1}{2}$,
\begin{eqnarray}\label{thm1-eq2bis}
n^{\nicefrac{(m+1)}{2}}\left[\widehat{X}^{(n)}_1-X^x_1\right]
\,{\stackrel{{\rm Law}}{\longrightarrow}}\,\,\,
\sigma(X^x_1)\left(
\int_0^1 h_{m}(X_s^x)[\sqrt{\mu_{2m+4}}\,dW_s+\mu_{m+3}\,dB_s]\right.\\
\left.+
\mu_{m+3}\int_0^1 g_m(X_s^x)ds\right),\nonumber
\end{eqnarray}

\quad with $W$ a Brownian motion independent of $B$. 

\textbullet\quad when $m$ is odd and  $H\in(\nicefrac{1}{2},1)$,
\begin{eqnarray}\label{314}
n^{(m+1)H}\left[\widehat{X}^{(n)}_1-X^x_1\right]
\,{\stackrel{{\rm Prob}}{\longrightarrow}}\,\,\,
\mu_{m+3}\,\sigma(X^x_1)\int_0^{B_1}h_m(\phi(x,y))dy.
\end{eqnarray}
\end{thm}
\begin{rem}\label{rm34}
{\rm 
\begin{enumerate}
\item For $m=0$ and $H>\nicefrac{1}{2}$, one recovers the convergence (\ref{euler-nn}). 
\item With the same method used to obtain (\ref{thm1-eq2bis}), one could prove (\ref{nn3})
with the help of Lemma 3.4 in \cite{N06}. Details are left
to the reader.
\item Actually, we could prove that the convergence is almost sure
 in (\ref{thm1-eq1}).
Also the convergence in (\ref{thm1-eq2}) is certainly almost sure, but the method of proof we have
used here does not allow to deduce it. Thus it remains an open question.
\item According to Theorem \ref{main-result}, whenever
$H\in(\nicefrac{1}{(m+2)},\nicefrac12)$ the scheme $\widehat{X}$ of size $m=2\kappa-1$ has the same
rate of convergence than the scheme $\widehat{X}$ of size $m=2\kappa$, namely $n^{(2\kappa+2)H-1}$. Thus,
it is a priori better to use odd-size schemes.
\item With the same method (see also Theorems 2 and 4 in \cite{NN}), one could also derive the exact rate
of convergence for the global error on the whole interval $[0,1]$. 
\end{enumerate}
}
\end{rem}

Observe that, under $(\mathscr{E})$,  the convergences (\ref{thm1-eq1}) 
and (\ref{thm1-eq2}) give the right lower error bound if
the probability that the right-hand side vanishes is strictly less than 1.
Due to (\ref{X}) and the fact that $B_t$ has a Gaussian density for any $t\in ]0,1]$, it is easy
to see that this last fact is equivalent to have that the real function
inside the integral, say $f_m$, is not identicaly zero.
Indeed, 
if $\int_0^1 g(B_s)ds=0$ almost surely for a certain $g\in{\rm C}^1_b(\mathbb{R})$, then
$0=D_u \int_0^1 g(B_s)ds = \int_u^1 g'(B_s)ds$, for any $u\in [0,1]$
(here $D$ denotes the Malliavin derivative with respect to $B$). We deduce that $g'(B_u)=0$, 
for any $u\in [0,1]$, and, since the support of the law of $B_1$ (for instance) is
$\mathbb{R}$, 
we obtain $g'= 0$. The desired conclusion follows easily. 

Except for $m=0$, solving $f_m= 0$ seems complicated. Nevertheless, when $m=1$, 
we can state:
\begin{prop}\label{prop1}
Assume that $(\mathscr{E})$ is in order, and moreover that $\sigma$ does not vanish.
Then the function $3\sigma'^3+6\sigma\sigma'\sigma''+\sigma^2\sigma'''$
(which is, up to a constant, the function appearing inside the integral
of the right-hand side of (\ref{thm1-eq2}) when $m=1$) 
is not identically zero.
\end{prop}

\begin{rem}
{\rm 
\begin{enumerate}
\item When $\sigma(x)=\sigma$ is constant, we have $\widehat{X}^{(n)}_1=X_1^x=x+\sigma\,B_1$. 
Consequently, the study of the rate of convergence in the case where $\sigma$ is a constant function 
is not interesting.
\item A corollary of Theorem \ref{main-result} and Proposition \ref{prop1} is that, 
under the additional hypothesis that $\sigma'$ does not vanish, 
the upper and lower error bounds always come from a convergence in probability
whenever $H\neq \nicefrac12$ and $m=1$. In particular, we never observe a phenomenon of the
type (\ref{nn2}). 
\end{enumerate}
}
\end{rem}

\noindent
{\bf Proof of Proposition \ref{prop1}}.
Since $\sigma'$ does not vanish, we have either $\sigma'>0$ or 
$\sigma'<0$. Suppose for instance that $\sigma'>0$, the proof for the other situation
being similar. Assume for a moment that $f:=3\sigma'^3+6\sigma\sigma'\sigma''+\sigma^2\sigma'''$ is identicaly zero. We then have $3\sigma'(\sigma^2)''=-\sigma(\sigma^2)'''$. We deduce that the derivative of 
$\sigma^3(\sigma^2)''$ is zero 
and then $(\sigma^2)''=\alpha\,\sigma^{-3}$ on $\R$, for some
$\alpha\neq 0$. Set $h=\sigma^2$; we have $h''\,h'=\alpha\,h'\,h^{-\nicefrac{3}{2}}$ or,
equivalently, $h'^2=-4\alpha\,h^{-\nicefrac{1}{2}}+\beta$ for some $\beta\in\mathbb{R}$. In particular,
we have $\beta-4\alpha\,y^{-\nicefrac{1}{2}}>0$, for any $y\in h(\R)$. Let $F$ be defined on $h(\R)$
by
$$
F(y)=\int_{h(0)}^y \frac{dz}{\sqrt{\beta-4\alpha\,z^{-\nicefrac{1}{2}}}}.
$$
For all $x\in \R$, we have
\begin{equation}\label{fhx}
F(\sigma(x)^2)=F(h(x))=x+\gamma,\mbox{ for some $\gamma\in\mathbb{R}$}.
\end{equation}
The function $\sigma^2$ being bounded, we necessarily have 
$h(x)\rightarrow\left(\nicefrac{4\alpha}{\beta}\right)^2$
(in particular $\beta\neq 0$)
as $x\to\infty$. Then, since $h''=\alpha\,h^{-\nicefrac32}$, this implies that $h''(x)\to\nicefrac{\beta^3}{(4^{3}\alpha^2)}$ as $x\rightarrow\infty$, 
which is in contradiction
with the fact that $h=\sigma^2$ is bounded. The proof of the proposition is done.\hfill$\Box$

\subsection{Proof of Theorem \ref{main-result}}

Here, and for the rest of the paper, we assume that $H$ belongs to $(\nicefrac{1}{(m+2)},1)$ and we denote
$\Delta_n = \max_{k=0, \ldots, n-1} |\Delta B_{\nicefrac{k}{n}}|$. 
We split the proof of Theorem \ref{main-result} into several steps.
\bigskip

{\it 1. General computations}.
The following lemma can be shown by using the same method
as in the proof of Lemma \ref{lm1}, but with Lebesgue integral instead of 
Newton-Cotes integral (and by taking into account that $\sigma\in{\rm C}^{\infty}(\mathbb{R})$
is bounded with bounded derivatives, in order to have uniform estimates):
\begin{lem}\label{lm2}
As $y\to 0$, we have, uniformly in $x\in\R$,
$$
\phi(x,y)=x+\sum_{j=0}^{m+2}\,\frac{1}{(j+1)!}\,\mathscr{D}^j\sigma(x)\,y^{j+1}+O(y^{m+4}).
$$
\end{lem}

\noindent
By applying this lemma to $x=\widehat{X}^{(n)}_{\nicefrac{k}{n}}$
and $y=\Delta B_{\nicefrac{k}{n}}$, we obtain, using the definition of 
$\widehat{X}^{(n)}_{\nicefrac{(k+1)}{n}}$,
\begin{multline*}
\widehat{X}^{(n)}_{\nicefrac{(k+1)}{n}}
=
\phi(\widehat{X}^{(n)}_{\nicefrac{k}{n}},
\Delta B_{\nicefrac{k}{n}})
-
\mathscr{D}^{m+1}\sigma
(\widehat{X}^{(n)}_{\nicefrac{k}{n}})\frac{\big(\Delta B_{\nicefrac{k}{n}}\big)^{m+2}}{(m+2)!}\\
-\mathscr{D}^{m+2}\sigma(\widehat{X}^{(n)}_{\nicefrac{k}{n}})\frac{\big(\Delta B_{\nicefrac{k}{n}})^{m+3}}{(m+3)!}
+O\big(\Delta_n^{m+4}\big).
\end{multline*}
By straightforward computations we get\footnote{
In fact, we rather obtain 
\begin{equation*}
\widehat{X}^{(n)}_{\nicefrac{(k+1)}{n}}
=
\phi\left(\widehat{X}^{(n)}_{\nicefrac{k}{n}},\Delta B_{\nicefrac{k}{n}}
+h_m(\widehat{X}^{(n)}_{\nicefrac{k}{n}})\big(\Delta B_{\nicefrac{k}{n}}\big)^{m+2}
+g_m(\widehat{X}^{(n)}_{\nicefrac{k}{n}})\big(\Delta B_{\nicefrac{k}{n}}\big)^{m+3}\right)
+O
\big(\Delta_n^{m+4}\big),
\end{equation*}
which is not exactly (\ref{diabolique}). But, in order to introduce 
$O\big(\Delta_n^{m+4}(B)\big)$ in the argument of $\phi$, we proceed as follows, by 
using the ellipticity property in hypothesis ($\mathscr{E}$):
$$
\phi(x,z)+O(\delta)=\phi(x,\phi^{-1}(x,\phi(x,z)+O(\delta)))=\phi(x,z+O(\delta)).
$$
}
\begin{equation}\label{diabolique}
\widehat{X}^{(n)}_{\nicefrac{(k+1)}{n}}
=
\phi\left(\widehat{X}^{(n)}_{\nicefrac{k}{n}},\Delta B_{\nicefrac{k}{n}}
+h_m(\widehat{X}^{(n)}_{\nicefrac{k}{n}})\big(\Delta B_{\nicefrac{k}{n}}\big)^{m+2}
+g_m(\widehat{X}^{(n)}_{\nicefrac{k}{n}})\big(\Delta B_{\nicefrac{k}{n}}\big)^{m+3}
+O
\big(\Delta_n^{m+4}\big)
\right),\vspace{0.5cm}
\end{equation}
with $g_m$ and $h_m$ given by (\ref{gm}).
By applying the group property (\ref{sgp}) repeatedly, we finally obtain that, for any $\ell\in\{1,\ldots,n\}$:
\begin{equation}\label{use-flow}
\widehat{X}^{(n)}_{\nicefrac{\ell}{n}}
=
\phi\left(x,B_{\nicefrac{\ell}{n}}
+\sum_{k=0}^{\ell-1}h_m(\widehat{X}^{(n)}_{\nicefrac{k}{n}})\big(\Delta B_{\nicefrac{k}{n}}\big)^{m+2}
+\sum_{k=0}^{\ell-1}g_m(\widehat{X}^{(n)}_{\nicefrac{k}{n}})\big(\Delta B_{\nicefrac{k}{n}}\big)^{m+3}
+O
\big(n\Delta_n^{m+4}\big)
\right).
\end{equation}
Since $\nicefrac{\partial\phi}{\partial y}=\sigma\circ\phi$ is bounded and $(\mathscr{E})$ is in order,
we deduce, as $n\to\infty$, 
\begin{equation}\label{wearat}
\sup_{\ell\in\{1,\ldots,n\}}\left|\widehat{X}_{\nicefrac{\ell}{n}}^{(n)}-X^x_{\nicefrac{\ell}{n}}\right|
=\sup_{\ell\in\{1,\ldots,n\}}\left|\widehat{X}_{\nicefrac{\ell}{n}}^{(n)}-\phi(x,B_{\nicefrac{\ell}{n}})\right|
=O
\big(n\Delta_n^{m+2}\big).
\end{equation}
In particular, $\widehat{X}_1^{(n)}$ converges almost surely to $X_1^x$ as $n\to\infty$,
since $H>\nicefrac{1}{(m+2)}$.
\bigskip

{\it 2. Proof of (\ref{thm1-eq1}).}
Let $m$ be an even integer.
As a consequence of (\ref{use-flow}) and (\ref{wearat}), 
we can write
\begin{equation}\label{yenamarredelavion}
\widehat{X}^{(n)}_{1}
=
\phi\left(x,B_{1}
+\sum_{k=0}^{n-1}h_m(X^{x}_{\nicefrac{k}{n}})\big(\Delta B_{\nicefrac{k}{n}}\big)^{m+2}
+O\big(n\Delta_n^{m+3}\big)+O\big(n^2\Delta_n^{2m+4}\big)
\right).
\end{equation}
(Observe that the main difference between (\ref{use-flow}) and the previous identity
is that the argument of $h_m$ is here $X^{x}_{\nicefrac{k}{n}}$).
By using (\ref{2-pow}) with $\kappa=m+2$, and due to the fact that $X_t^x=\phi(x,B_t)$
and $\nicefrac{\partial \phi}{\partial y}=\sigma\circ\phi$, we finally obtain 
(\ref{thm1-eq1}).
\bigskip

{\it 3. Proof of (\ref{thm1-eq2}) for $H\in(\nicefrac{2}{(2m+3)},\nicefrac{1}{2})$.}
Let $m$ be an odd integer and assume that $H\in(\nicefrac{2}{(2m+3)},\nicefrac{1}{2})$.
Thanks to (\ref{wearat}), identity
(\ref{use-flow}) can be transformed into
\begin{multline}\label{eq78}
\widehat{X}^{(n)}_{\nicefrac{\ell}{n}}
=
\phi\left(x,B_{\nicefrac{\ell}{n}}
+\sum_{k=0}^{\ell-1}h_m(\widehat{X}^{(n)}_{\nicefrac{k}{n}})\big(\Delta B_{\nicefrac{k}{n}}\big)^{m+2}\right.\\
\left.+\sum_{k=0}^{\ell-1}g_m({X}^{x}_{\nicefrac{k}{n}})\big(\Delta B_{\nicefrac{k}{n}}\big)^{m+3}
+O\big(n\Delta_n^{m+4}\big)+O\big(n^2\Delta_n^{2m+5}\big)\right).
\end{multline}
On the other hand, due to $(\mathscr{E})$, we have, 
for any fixed $M\geq 1$ and uniformly in $x\in\R$,
\begin{equation*}
 \phi(x,y_{2})= \phi(x,y_{1}) +  \sum_{j=1}^{M}
\frac{1}{j!} \frac{ \partial^{j} \phi}{ \partial y^{j}} 
(x,y_{1}) (y_{2}-y_{1})^{j} + O((y_{2}-y_{1})^{M+1}  ). 
\end{equation*}
Combined with (\ref{use-flow}), it yields
\begin{multline}\label{eq_re_in_1}
{\widehat X}^{(n)}_{k/n}=
X_{k/n}^x+ \sum_{j=1}^{M} 
\frac{1}{j!} \frac{\partial^{j} \phi}{\partial y^{j}}(x, B_{k/n})
\left( \sum_{k_1=0}^{k-1}h_m(\widehat{X}^{(n)}_{k_1/n})
\big(\Delta B_{\nicefrac{k_1}{n}}\big)^{m+2} 
+O ( n \Delta_n^{m+3})
\right)^{j} \\ 
+ O(n^{M+1} \Delta_n^{(m+2)(M+1)}).
\end{multline} 
By using (\ref{eq_re_in_1}) with $M=1$ 
as well as the equality $\nicefrac{\partial \phi }{\partial y}=\sigma\circ\phi$,
we get
$$
{\widehat X}^{(n)}_{\nicefrac{k}{n}}=
X_{\nicefrac{k}{n}}^x
+\sigma(X^x_{\nicefrac{k}{n}})\sum_{k_1=0}^{k-1}h_m(\widehat{X}^{(n)}_{\nicefrac{k_1}{n}})
\big(\Delta B_{\nicefrac{k_1}{n}}\big)^{m+2}
+ 
O
\big(n^2\Delta_n^{2m+4}\big)
+O
\big(n\Delta_n^{m+3}\big)
$$
and then, by (\ref{wearat}),
$$
{\widehat X}^{(n)}_{\nicefrac{k}{n}}=
X_{\nicefrac{k}{n}}^x
+\sigma(X_{\nicefrac{k}{n}}^x)\sum_{k_1=0}^{k-1}h_m(X^x_{\nicefrac{k_1}{n}})
\big(\Delta B_{\nicefrac{k_1}{n}}\big)^{m+2}
+ 
O
\big(n^2\Delta_n^{2m+4}\big)
+O
\big(n\Delta_n^{m+3}\big).
$$
By inserting the previous equality in (\ref{eq78})
with $\ell=n$, we obtain
\begin{multline}\label{eq34}
\widehat{X}^{(n)}_{1}
=
\phi(x,B_{1}
+\sum_{k=0}^{n-1}h_m({X}^{x}_{\nicefrac{k}{n}})\big(\Delta B_{\nicefrac{k}{n}}\big)^{m+2}
+\sum_{k=0}^{n-1}g_m({X}^{x}_{\nicefrac{k}{n}})\big(\Delta B_{\nicefrac{k}{n}}\big)^{m+3}\\
+\sum_{k=0}^{n-1}\sigma h'_m({X}^{x}_{\nicefrac{k}{n}})\big(\Delta B_{\nicefrac{k}{n}}\big)^{m+2}
\sum_{j=0}^{k-1}h_m({X}^{x}_{\nicefrac{j}{n}})\big(\Delta B_{\nicefrac{j}{n}}\big)^{m+2}\\
+O
\big(n^3\Delta_n^{3m+6}\big)+O\big(n^2\Delta_n^{2m+5}\big)
+O
\big(n\Delta_n^{m+4}\big)
).
\end{multline}
Due to (\ref{3-pow}) with $\kappa=m+2$ we have, as $n\to\infty$,
$$
n^{(m+3)H-1}\sum_{k=0}^{n-1}h_m(X^{x}_{\nicefrac{k}{n}})\big(\Delta B_{\nicefrac{k}{n}}\big)^{m+2}
\,\,\,{\stackrel{{\rm Prob}}{\longrightarrow}}\,\,\,
-\frac{\mu_{m+3}}{2}\,\int_0^1 \sigma h'_m(X^x_s)ds
$$
and also, due this time to (\ref{2-pow}) with $\kappa=m+3$, as $n\to\infty$,
$$
n^{(m+3)H-1}\sum_{k=0}^{n-1}g_m(X^{x}_{\nicefrac{k}{n}})\big(\Delta B_{\nicefrac{k}{n}}\big)^{m+3}
\,\,\,{\stackrel{{\rm Prob}}{\longrightarrow}}\,\,\,
\mu_{m+3}\,\int_0^1  g_m(X^x_s)ds.
$$
Moreover, since we assume in this step that $H>\nicefrac{2}{(2m+3)}$, we have, as $n\to\infty$,
$$
n^{(m+3)H}\Delta_n^{m+4}
\,\,\,{\stackrel{{\rm Prob}}{\longrightarrow}}\,\,\,0,\quad
n^{(m+3)H+1}\Delta_n^{2m+5}
\,\,\,{\stackrel{{\rm Prob}}{\longrightarrow}}\,\,\,0\quad\mbox{and}\quad
n^{(m+3)H+2}\Delta_n^{3m+6}
\,\,\,{\stackrel{{\rm Prob}}{\longrightarrow}}\,\,\,0.
$$
At this level, we need the following result which is contained in \cite[Proposition 7]{NNT}:
\begin{lem}
Fix an integer $q\geq 2$ and denote by $H_{q}$ the $q$th Hermite polynomial. 
Let $f\in C^{2q}(\mathbb{R})$ be bounded with bounded derivatives and, for 
$k\in\{1,\ldots,n\}$, denote
\begin{equation*}
S_k^{(q)}(f):=\sum_{j=0}^{k-1} f(B_{\nicefrac{j}{n}})H_q(n^H\Delta B_{\nicefrac{j}{n}}).
\end{equation*}
Then 
\begin{equation}\label{dollardollar}
E\big|S_k^{(q)}(f)\big|^2=O(n^{1\vee (2-2Hq)})\quad\mbox{as $n\to\infty$, uniformly in $k$}.\\
\end{equation}
\end{lem}
Recall also that, since $m+2$ is odd, the monomial $x^{m+2}$ may be expanded in terms of the Hermite polynomials as follows:
\begin{equation}\label{her-odd}
x^{m+2}-\mu_{m+3} x=\sum_{q=1}^{\nicefrac{(m+1)}2}a_{m+2,2q+1}H_{2q+1}(x),\quad
\mbox{for some universal constants $a_{m+2,2q+1}$}.
\end{equation}
Therefore,
for $k\in\{1,\ldots,n\}$,
\begin{multline*}
n^{(m+3)H-1}\sum_{j=0}^{k-1} h_m(X^{x}_{\nicefrac{j}{n}})\big(\Delta B_{\nicefrac{j}{n}}\big)^{m+2}
-\mu_{m+3}\,n^{2H-1}\sum_{j=0}^{k-1} 
h_m(X^{x}_{\nicefrac{j}{n}})\Delta B_{\nicefrac{j}{n}}\\
=n^{H-1}\sum_{q=1}^{\nicefrac{(m+1)}2}a_{m+2,2q+1}S_k^{(2q+1)}\big(h_m(\phi(x,\cdot))\big)
\end{multline*}
and, by (\ref{dollardollar}),
\begin{eqnarray*}
&&E\left|n^{(m+3)H-1}\sum_{j=0}^{k-1} h_m(X^{x}_{\nicefrac{j}{n}})\big(\Delta B_{\nicefrac{j}{n}}\big)^{m+2}
-\mu_{m+3}\,n^{2H-1}\sum_{j=0}^{k-1} 
h_m(X^{x}_{\nicefrac{j}{n}})\Delta B_{\nicefrac{j}{n}} \right|^2\\
&&\hskip8cm= 
O\big(n^{(2H-1)\vee(-4H)}\big).
\end{eqnarray*}
Hence
\begin{eqnarray*}\label{eq91}
&&E\left|n^{(m+3)H-1}\sum_{k=0}^{n-1}\sigma h'_m(X^x_{\nicefrac{k}{n}})
\big(\Delta B_{\nicefrac{k}{n}}\big)^{m+2}
\sum_{j=0}^{k-1} h_m(X^{x}_{\nicefrac{j}{n}})\big(\Delta B_{\nicefrac{j}{n}}\big)^{m+2}\right.\\
&&\hskip1cm \left. 
-\mu_{m+3}\,n^{2H-1}\sum_{k=0}^{n-1}\sigma h'_m(X^x_{\nicefrac{k}{n}})
\big(\Delta B_{\nicefrac{k}{n}}\big)^{m+2}
\sum_{j=0}^{k-1} h_m(X^{x}_{\nicefrac{j}{n}})\Delta B_{\nicefrac{j}{n}}
\right|\\
&&\hskip8cm
=O\big(n^{(\frac12-H-mH)\vee (1-4H-mH)}\big)
\end{eqnarray*}
which tends to zero
as $n\rightarrow\infty$,
because $H>\nicefrac{2}{(2m+3)}$ implying $H>\nicefrac{1}{(2m+2)}$ and
$H>\nicefrac{1}{(m+4)}$.
Moreover, by the mean value theorem:
$$
\sum_{j=0}^{k-1} h_m(X^{x}_{\nicefrac{j}{n}})\Delta B_{\nicefrac{j}{n}}
=\int_0^{B_{\nicefrac{k}n}} h_m\big(\phi(x,z)\big)dz-\frac12 \sum_{j=0}^{k-1} \sigma 
h'_m(X^{x}_{\theta_{\nicefrac{j}{n}}})\big(\Delta B_{\nicefrac{j}{n}}\big)^{2},
$$
for some
$\theta_{\nicefrac{j}{n}}$ between $\nicefrac{j}{n}$ and $\nicefrac{(j+1)}{n}$.
Consequently, since $H<1/2$, we have
$$
E\left|\sum_{j=0}^{k-1} h_m(X^{x}_{\nicefrac{j}{n}})\Delta B_{\nicefrac{j}{n}}
\right|^2=O(n^{2-4H})
$$
so that
$$
E\left|n^{2H-1}\sum_{k=0}^{n-1}\sigma h'_m({X}^{x}_{\nicefrac{k}{n}})\big(\Delta B_{\nicefrac{k}{n}}\big)^{m+2}
\sum_{j=0}^{k-1}h_m({X}^{x}_{\nicefrac{j}{n}})\Delta B_{\nicefrac{j}{n}}
\right|=O(n^{-(m+2)H})\longrightarrow 0.
$$
Finally, by combining all these convergences to zero together, we get
$$
n^{(m+3)H-1}\sum_{k=0}^{n-1}\sigma h'_m({X}^{x}_{\nicefrac{k}{n}})\big(\Delta B_{\nicefrac{k}{n}}\big)^{m+2}
\sum_{j=0}^{k-1}h_m({X}^{x}_{\nicefrac{j}{n}})\big(\Delta B_{\nicefrac{j}{n}}\big)^{m+2}\,\,\overset{{\rm Prob}}{\longrightarrow}\,\,0,
$$
so that
the proof of (\ref{thm1-eq2}) is done in the case where
$H>\nicefrac{2}{(2m+3)}$.
\bigskip

{\it 4. Proof of (\ref{thm1-eq2}) for $H\in(\nicefrac{1}{(m+2)},\nicefrac{2}{(2m+3)}]$.}
It suffices to use (\ref{eq_re_in_1}) with the appropriate $M$ for the considered $H$ and then
to proceed as in the previous step. The remaining details are left to the reader.
\bigskip

{\it 5. Proof of (\ref{thm1-eq2bis}).} By going one step further in 
(\ref{eq78}) using (\ref{wearat}), we get
\begin{multline*}
\widehat{X}^{(n)}_{1}
=
\phi(x,B_{1}
+\sum_{k=0}^{n-1}h_m(X^{x}_{\nicefrac{k}{n}})\big(\Delta B_{\nicefrac{k}{n}}\big)^{m+2}
+\sum_{k=0}^{n-1}g_m(X^{x}_{\nicefrac{k}{n}})\big(\Delta B_{\nicefrac{k}{n}}\big)^{m+3}\\
+
O\big(n^2\Delta_n^{2m+4}\big)
+O\big(n\Delta_n^{m+4}\big)).
\end{multline*}
Whenever $m\ge 3$ and since $H=\nicefrac{1}{2}$, we have, as $n\to\infty$:
$$
n^{\nicefrac{(m+1)}{2}+1}\Delta_n^{m+4}
\,\,\,{\stackrel{{\rm Prob}}{\longrightarrow}}\,\,\,0
\quad\mbox{and}\quad
n^{\nicefrac{(m+1)}{2}+2}\Delta_n^{2m+4}
\,\,\,{\stackrel{{\rm Prob}}{\longrightarrow}}\,\,\,0.
$$
Hence, for $m\ge 3$, (\ref{thm1-eq2bis}) is an immediate consequence of 
(\ref{mb3-pow}) and of the previous two relations.

If $m=1$, we rather need to use (\ref{eq34}). Since $H=\nicefrac{1}{2}$,
we have, as $n\to\infty$:
$$
n^{2}\Delta_n^{5}
\,\,\,{\stackrel{{\rm Prob}}{\longrightarrow}}\,\,\,0,
\quad
n^{3}\Delta_n^{7}
\,\,\,{\stackrel{{\rm Prob}}{\longrightarrow}}\,\,\,0
\quad\mbox{and}\quad
n^{4}\Delta_n^{9}
\,\,\,{\stackrel{{\rm Prob}}{\longrightarrow}}\,\,\,0.
$$
Finally, combining these convergences with (\ref{2-pow}) (for $H=\nicefrac{1}{2}$), (\ref{mb3-pow}) and 
the fact that
\begin{multline*}
{\rm E}\left|
n\sum_{k=0}^{n-1}\sigma h'_1(X^x_{\nicefrac{k}{n}})
\big(\Delta B_{\nicefrac{k}{n}}\big)^{3}
\sum_{j=0}^{k-1} 
h_1(X^{x}_{\nicefrac{j}{n}})\big(\Delta B_{\nicefrac{j}{n}}\big)^{3}
\right|^2\\
=n^2\sum_{k=0}^{n-1}\sum_{j=0}^{k-1}{\rm E}\left|
\sigma h'_1(X^x_{\nicefrac{k}{n}})
\big(\Delta B_{\nicefrac{k}{n}}\big)^{3}
h_1(X^{x}_{\nicefrac{j}{n}})\big(\Delta B_{\nicefrac{j}{n}}\big)^{3}
\right|^2
=O(n^{-2})\longrightarrow 0\quad\mbox{as $n\rightarrow\infty$},
\end{multline*}
we obtain (\ref{thm1-eq2bis}) also for $m=1$.
\bigskip

{\it 6. Proof of (\ref{314})}. By combining (\ref{yenamarredelavion}) with the fact that
$X_t^x=\phi(x,B_t)$ and $\nicefrac{\partial\phi}{\partial y}=\sigma\circ\phi$, we get
$$
\widehat{X}_1^{(n)}=X_1^x +\sigma(X_1^x)\sum_{k=0}^{n-1}
h_m(X^x_{\nicefrac{k}{n}})
\big(\Delta B_{\nicefrac{k}{n}}\big)^{m+2}
+O\big(n\Delta_n^{m+3}\big)+O(n^2\Delta_n^{2m+4}\big).
$$
Since $H>\nicefrac12$, we have, as $n\to\infty$,
$$
n^{(m+1)H+1}\Delta_n^{m+3}
\,\,\,{\stackrel{{\rm Prob}}{\longrightarrow}}\,\,\,0
\quad\mbox{and}\quad
n^{(m+1)H+2}\Delta_n^{2m+4}
\,\,\,{\stackrel{{\rm Prob}}{\longrightarrow}}\,\,\,0.
$$
Hence (\ref{314}) is an immediate consequence of (\ref{28bis}).
\fin

\noindent
{\bf Acknowledgement}: We thank Jean Jacod for helpful discussion about 
stable convergence and especially about reference \cite{jacod}
and Serge Cohen for its suggestion of reference \cite{KG}.

\end{document}